\documentclass[a4paper,12pt, leqno]{article}
\usepackage{amsmath, amssymb}
\usepackage{txfonts}
\begin{document}
\thispagestyle{empty}
\begin{center}
{\it\huge
On the classical solutions of the three dimensional Stefan problem
}
\end{center}
\newpage
\thispagestyle{empty}
\mbox{}
\newpage
\setcounter{page}{1}
\begin{center}
{\it\Large
On the classical solutions of the three dimensional Stefan problem
}
\end{center}
\begin{center}
{\large
Jun-ichi Koga
}
\end{center}
\begin{center}
{\bf
1.Introduction
}
\end{center}
\quad We study the Stefan problem which corresponds to the interface  between the ice and water. The problem is following. {\it Do the classical solutions in time exists on a boundary especially for the three dimensional Stefan problem?} Because our difficulty of analysis is the jump-behavior on a boundary. We investigate the following story of the global existence theorem of the Stefan problem.\\
\quad In physics, the history of the Stefan problem goes back to Lam\'e and Clapeyron [7]. They proposed the fundamental properties for the Stefan problem in nature. Furthermore Friedman said in [3] that much harder is the problem in the case of several space variables. Kamenomostskaya slightly proved the classical Stefan problem in [3] . In this paper the author proves the existence and uniqueness of a generalized solution for the three-dimensional Stefan problem, later improved by her master Olga Oleinik. Salsa [9], on the other hand, suggested recent results and open problems of the Stefan problem. According to [9], \\
$$
\mathcal{L} u=\mathrm{Tr} (A(x,t)D^2u)+b(x,t)\cdot\nabla u
$$
are difficult to consider with the Stefan conditions.\\
\quad One of the famous open Stefan problems is to prove that Lipschitz (but non necessarily flat in space) free
boundaries are smooth, under a nondegeneracy condition. This result would allow, for instance, to treat nonlinear divergence operators of the type
$$
u_t-\mathrm{div} (A(x,t,u))
$$
with Lipschitz continuous coefficients. The main difficulty here stems form the construction of the continuous family of deformations constructed to decrease $\varepsilon$ in the $\varepsilon$-monotonicity conditions. In fact, in this construction the flatness of the spatial sections of $F(u)$ plays a major role. Reutskiy [8], however, calculate the Stefan problems numerically with the method of approximate fundamental solutions (MAFS). \\
\quad Then, due to Friedman [3,4,6], set $\Omega_T=G\times (0,T)$ and we introduce the elliptic operator
\begin{equation}
L_i=\sum_{j,k=1}^n a_{jk}^i (x,t)\frac{\partial^2}{\partial x_j\partial x_k}+\sum_{j=1}^nb^i_j (x,t)\frac{\partial}{\partial x_j}+c^i(x,t),\tag{$1.1$}
\end{equation} 
where coefficients satisfying: $a^i_{jk}$, $\nabla_xa^i_{jk}$, $\nabla^2_xa^i_{jk}$, $b^i_j$, $\nabla_xb^i_j$, and $c^i$ are continuous in $\overline{\Omega}_\infty$.\\
\quad Now, we consider the Stefan problem in several space variables. More precisely, the system of the parabolic equations are the following equations using $(1.1)$ in the sense of Petrovskii
\begin{equation}
\left\{
\begin{array}{lcl}
L_i u=u_t & & 0<x<s(t),\,\,t>0\\\\
u(0,t)=f(t)     & & f(t)\ge 0\,\, t>0\\\\
u(\vec{s}(t),t)=0 & & t>0\,\,\vec{s}(0)=b\\\\
v=\left(\nabla u_S-\nabla u_L\right)\cdot\widehat{\mathbf{n}} & & t>0\cdots\cdots (*)
\end{array}
\right.\tag{$1.2$}
\end{equation}
however, ``$t\geq 0$'' for $t$ holds by Lemma 2.1.\\\\
{\bf
Definition 1.1
}\quad We call the {\it Fridrich's mollifier} 
$$
f^\epsilon (x)=\int_U\eta_\epsilon (x-y)f(y)dy=\int_{B(0,\epsilon )}\eta_\epsilon (y)f(x-y)dy
$$
for $x\in U_\epsilon$, where $\eta_\epsilon (x)$ is $\displaystyle\frac{1}{\epsilon^n}\eta\left(\frac{x}{\epsilon}\right)$ for $\eta\in C^\infty (\mathbb{R}^n)$.\\\\  
\quad The ``global'' (classical) solutions in time very surprisingly exists, proved by E. Hanzawa. However, we could get the classical solutions in the weak sense, namely $j_\epsilon$ is a Friedrich's mollifier and we introduce this operator the convolution operator as $j_\epsilon *$ like the elliptic operators, to prove that the function $j_\epsilon*f$ ($f\in L^2$) is $C_0^\infty$. \\\\
{\bf
Theorem 1.2
}\quad Assume that $j_\epsilon*$ is a Friedrich's mollifier. Provided that $u$ is a unique solution of $(1.2)$, and that $u\in L^2(\Omega)$, $j_\epsilon *u\in C_0^\infty$.\\\\
{\bf 
Remark 1.3
}
\quad Let $\displaystyle V_n:=k[(\partial\rho / \partial n)]_\Gamma\,\,(\mathrm{for\ any\ }\rho\in C_0^\infty)$ be the Stefan condition with jump. This is one of the most difficult points in our problem. So we must calculate this equality with the velocities of the interface to obtain our proof of the classical solutions.\\\\
{\bf
Remark 1.4
}
\quad For one dimensional Stefan problem, the Stefan condition $(*)$ is smooth in [4], but for $n\geq 2$, this condition $(*)$ is jump on the boundary.\\\\  
\quad The layout of our paper is as follows: in Chapter 1, we give Introduction. In Chapter 2, we state and prove preliminaries, Chapter 3 is devoted to the proof of our main theorem, and finally Chapter 4 contains Appendix.
\begin{center}
{\bf
2.Preliminaries
}
\end{center}
\quad Before starting to state and to prove our lemmas, we introduce the special parabolic type (which we call);
\begin{equation}
v^m(x,t)=u(x,t)-\frac{1}{8n}\left(\left|x-x^m\right|^2+\left(t^m-t\right)\right),\tag{$2.1$}
\end{equation}
due to [1], where $u(x,t)$ is a solution of  a normal parabolic heat equation;
\begin{equation}
d\left(\left(x,t\right),\left(x_0,t_0\right)\right)=\left(\left|x-x_0\right|^2+\left|t-t_0\right|\right)^{1/2}.\tag{$2.2$}
\end{equation}
{\bf Lemma.2.1}\quad $u_t$ is continuous at $t=0$.\\\\
{\bf Proof.}\quad Define
\begin{equation}
G(t)=\left\{x\in\mathbb{R}^n\left\backslash G_0\left|\right.\right. u(x,t)>0\right\}\tag{$2.3$}
\end{equation}
and the more detail part of our proof based on [3] was claimed that 
\begin{equation}
G(t)\,\, \mathrm{is\,\, contained\,\, in\,\, a}\,\,\delta(t)\mathrm{-neighborhood\,\, of}\,\, G,\tag{$2.4$}
\end{equation}
where $\delta (t)\to 0$ as $t\to 0$. In the set $\mathbf{R}$ of points $(x,t)$ such that $u(x,t)>0$ and
\begin{equation*}
\left|x-x^m\right|<c,\,\, 0<t<t^m.
\end{equation*}
Clearly, by the simple calculations, we get
\begin{equation*}
\begin{array}{rcl}
v^m(x^m,t^m) &=&\displaystyle u(x^m,t^m)-\frac{1}{8n}\left(\left|x^m-x^m\right|^2+\left(t^m-t^m\right)\right)\\\\
                   &=& u(x^m,t^m)\\\\
                   &>& 0,
\end{array}
\end{equation*}
notice that this equations are simple but we can get the validity of the strong maximum principle through this equality, because of the assumption and the following inequalities are obtained,
\vspace{20truecm}
\begin{equation*}
\begin{array}{rcl}
\Delta v^m(x,t) &=&\displaystyle \Delta u - \frac{1}{8n}\left(\Delta\left|x-x^m\right|+\Delta\left(t^m-t\right)\right)\\\\
                     &=& \displaystyle\Delta u-\frac{1}{16n}\\\\
v^m_t(x,t)        &=& \displaystyle u_t(x,t)-\frac{1}{8n}\left(\left|x-x^m\right|_t+\left(t^m-t\right)_t\right)\\\\
                     &=& \displaystyle u_t-\frac{1}{8n}\\\\
\displaystyle\left(\Delta u-u_t\right)-\left(\frac{1}{16n}-\frac{1}{8n}\right) &=&\displaystyle 0+\frac{1}{8n}-\frac{1}{16n}\\\\
                     &=& \displaystyle\frac{1}{16n}\\\\
                     &>& 0
\end{array}
\end{equation*}
in $\mathbf{R}$. So $\Delta v^m -v^m_t>0$ in $\mathbf{R}$. By the maximum principle, $v$ must take a 
positive maximum value on the parabolic boundary of $\mathbf{R}$, that is to say, at $\left(x^m_1,t^m_1\right)$. Thus 
\begin{equation}
u\left(x^m_1,t^m_1\right)>\frac{1}{8n}\left(\left|x^m-x^m_1\right|^2+\left(t^m-t\right)\right)\tag{$2.5$}
\end{equation}
Since $t_1^m=0$ is impossible and since $(x_1^m,t_1^m)$ can not lie on the free boundary, we must have $\left|x^m-x^m_1\right|=c$. Hence (2.5) gives
$$
u(x^m_1,t^m_1)>\frac{c^2}{8n},\,\, t_1^m\to 0,\,\, \mathrm{dist}(x^m_1, G)>c
$$
which is impossible again since $h(x)=0$ if $x\not\in G$.\\
\quad Since $u_t\ge 0$
\begin{equation}
u(x,t)>0\,\,\,\, if\,\, x\in G.\tag{$2.6$}
\end{equation}
From (2.4), (2.6) and the boundedness of $u_t$ we deduce that
\begin{equation}
\int_{\mathbf{R}^n\backslash G_0}\left|u_t(x,t)-h(x)\right|dx\to 0\,\, if\,\, t\to 0;\tag{$2.7$}
\end{equation}
here we use the condition that $\Gamma$ is Lipschitz.\\
\quad We now proceed to prove the continuity of $u_t$ at a point $(y,0)$; it suffices to take y in $\Gamma$. Let $\mathbf{K}$ be a small ball with center $y$, and $w^\epsilon\ge 0$ be the solution of 
$$
\Delta w^\epsilon-w_t^\epsilon=0\,\, in\,\, \mathbf{K}\times (\epsilon , 1)
$$
$$
w^\epsilon=u_t\,\, \mathrm{on\,\, the\,\, parabolic\,\, boundary}\,\, \mathbf{K}\times(\epsilon , 1).
$$
By (2.7)
\begin{equation}
w^\epsilon (x,t)\to w^0(x,t) \,\,\mathrm{if}\,\, \epsilon\to 0.\tag{$2.8$} 
\end{equation}
\quad From Lemma 2.2 we have, for any $\epsilon>0$,
$$
u_t(x,t)\le w^\epsilon (x,t)
$$
Taking $\epsilon\to 0$ and using (2.8), we get
\begin{equation}
u_t(x,t)\le w^0(x,t)\tag{$2.9$}.
\end{equation} 
Since $w^0(x,0)=h(x)$ is continuous at $y$, $w^0(x,t)\to h(y)=0$ as $x\to y$, $t\to 0$.\,\, This completes the proof of Lemma 2.1. $\Box$\\\\
{\bf
Lemma 2.2
}\quad Let $w$ be a bounded measurable function in a cylinder $D\times (0,T)$, which is subcaloric (that is to say, $\Delta w-w_t\geq 0$ in the sense of the distribution). Then there exists a function $\widetilde{W}$ such that
\begin{itemize}
\item $\widetilde{w}=w$ a.e. in $d\times (0,T)$,
\item $\widetilde{w}$ is upper semicontinuous in $D\times (0,T)$,
\item for any ball $K, \overline{K}\subset D$, if $z$ satisfies
$$
\begin{array}{l}
\Delta z-z_t=0\,\,\mathrm{in\ } K\times (t_0,t_1)\,\, (0<t_0<t_1<T)\\\\
z=\widetilde{w}\,\,\mathrm{on\ the\ boundary\ of\ } K\times (t_0,t_1)
\end{array}
$$
\end{itemize}
then $z\leq \widetilde{w}$ in $K\times (t_0,t_1)$.
\begin{center}
{\bf
3. Outline of proof of Lemma 2.2
}
\end{center}
Let $C_\rho =(x^0,t^0)=\{(x,t);|x-x^0|<\rho, t^0-\rho^2<t<t^0\}$, $B_\rho(x^0)=\{x;|x-x^0|<\rho\}$ and denote by $G_\rho (x,t, \xi)$ the Green function for the heat operator with singularity at $x=\xi$, $t=0$. If $w$ is a smooth function then
$$
\begin{array}{rcl}
w(x^0,t^0) & = & \displaystyle\int_{B_\rho (x^0)}G_\rho(x^0,\rho^2,\xi)w(\xi, t^0-\rho^2)d\xi\\\\
              & = & \displaystyle -\int_{t_0-\rho^2}^{t^0}\int_{\partial B_\rho (x^0)}\frac{\partial G_\rho (x^0,\tau ,\xi)}{\partial n}w(\xi ,\tau)dS_\xi d\tau\\\\
             & = & \displaystyle -\int_{C_\rho (x^0,t^0)}G_\rho (x_0,\tau,\xi)\left(\Delta w-w_t\right)(\xi, \tau)d\xi d\tau.
\end{array}
$$
Integrating over $\rho$, $R<\rho <2R$, we get
$$
\begin{array}{rcl}
w(x^0,t^0) & = &\displaystyle w_R(x^0,t^0)-\frac{1}{R}\int_R^{2R}\int_{C_\rho (x^0,t^0)}G_\rho (x_0,\tau,\xi)\left(\Delta w-w_t\right)(\xi, \tau)d\xi d\tau\\
\end{array}
$$
where
$$
\begin{array}{rcl}
w_R(x^0,t^0) & = & \displaystyle\frac 1R\int^{2R}_R\int_{C_\rho (x_0,t_0)}G_\rho (x_0,\rho^2, \xi)w(\Delta w-w_\tau)(\xi, \tau)d\xi d\tau d\rho\\\\%\vspace{-5truemm}
& & \displaystyle -\frac{1}{R}\int_R^{2R}\int_{t^0-\rho^2}^{t^0}\int_{\partial B_\rho (x^0)}\frac{\partial G_\rho (x^0,\tau ,\xi)}{\partial n}w(\xi, \tau)dS_\xi d\tau d\rho.
\end{array}
$$
Now let $w$ be a bounded measurable function as in lemma, and define $w_R$ by the above equations. Notice that $w_R$ is continuous function. We claim that 
$$
w_R\leq w_{R^\prime}\,\,\mathrm{if\ }R<R^\prime
$$
To prove this, we use the Friedrich's mollifiers $j_\epsilon *$. Notice that any smooth function with these mollifiers has a certain support compact and that the operators we consider are subcaloric for $\rho,\rho^\prime >0\,\, (\rho<\rho^\prime)$,
$$
\left(j_\epsilon *w\right)_R\leq\left(j_\epsilon *w\right)_{R^\prime}\,\,\mathrm{if\ }R<R^\prime.
$$
Taking $\epsilon \to 0$, we obtain the desired claim.
\begin{center}
{\bf
4. Proof of Theorem 1.2
}
\end{center}
\quad The proof of Theorem 1 is similar to Yi's proof $[10]$ for the global classical solutions of the Muskat problem.\\ \quad {\it Step.1.}\quad For the Stefan condition, we simply compute
\begin{equation}
\begin{array}{ccccl}
V_n & = & \displaystyle -k_1\frac{\partial u}{\partial n} & = & -k_1\nabla u\cdot \mathbf{n}\\\tag{$3.1$}
& & & = & \displaystyle -k_1\left(\frac{\partial u}{\partial x}, \frac{\partial u}{\partial y}, \frac{\partial u}{\partial z}\right)\cdot\frac{1}{\sqrt{1+\rho_x^2+\rho_y^2}}\left(-\rho_x,-\rho_y, 1\right)\\
& & & = & \displaystyle \frac{\partial_t\rho}{\sqrt{1+\rho_x^2+\rho_y^2}}.
\end{array}
\end{equation}
Thus we get
\begin{equation}
\partial_t\rho =-k_1\left(u_z-\rho_xu_x-\rho_yu_y\right),\tag{$3.2$}
\end{equation}
which, therefore, rationalizes the heat equations with our conditions.\\
\quad {\it Step.2.}\quad More precisely, the differential scheme for the heat operator,
$$
\left\{
\begin{array}{rcl}
\displaystyle \frac{u^{(n)}_{i+1}-u^{(n)}_i}{\Delta t}-\frac{u^{(n)}_{i+1}+2u^{(n)}_i-u^{(n)}_{i-1}}{2(\Delta x)^2} & = & 0\\\\
\displaystyle u(0,t) & = &f(t)\\\\
\displaystyle u(\vec{s}(t),t) & = & 0\\\\
\displaystyle \partial_t\rho & = &-k_1\left(u_z-\rho_xu_x-\rho_yu_y\right)
\end{array}
\right.
$$
in H\"older spaces $C^{2,2}(\Omega)$. Then using  Navier-Stokes method in computational fluids dynamics (CFD) for this equations, we get the desired results. 
\\
\quad {\it Step.3.}\quad Lemma 2.1and the numerical analysis for the integral yields
\begin{equation}
\begin{array}{rcl}
\displaystyle\int_\Omega\int_0^t u_t-\Delta u dxdt& = &\displaystyle\int_\Omega\int_0^t\frac{\partial}{\partial t} udxdt-\int_\Omega\int_0^t\Delta udxdt,\tag{$3.3$}\\\\
& = & \displaystyle \int_\Omega\frac{\partial}{\partial t}\int_0^t udxdt-\int_\Omega\int_0^t\Delta u dxdt,\\\\
& = & \displaystyle \int_\Omega udx-\int_\Omega\int_0^t\Delta udxdt,\\\\
& = & 0
\end{array}
\end{equation}
and $u$ is, for the initial data $\phi$,
\begin{equation}
u(x,t)=\frac{1}{\sqrt{4\pi t}}\int_\Omega\phi(\xi)\exp\left(-\frac{|x-\xi|^2}{4t}\right)d\xi\,\, t\in\mathbb{R},x\in\mathbb{R}^n\tag{$3.4$}
\end{equation}
so we obtain the result $(3.3)$ of the infinitely continuously differentiable function $u$. From Step.2, we obtain the local existence in space. However we get the global one on this stage. This completes the proof of Theorem 2. 1.$\Box$\\\\
{\bf acknowledgments}\\\\
The author would like to thank $[1]$ for almost referred to the author's paper. 
\begin{center}
{\bf
References
}
\end{center}
\begin{itemize}
\item[{[1]}] Caffarelli, L., A., Friedman, A., Continuity of the temperature in the Stefan problem, Indiana Univ. Math. J. {\bf 28}, 1979, pp. 53-70.
\item[{[2]}] Caginalp, G, An analysis of a phase field model of a free boundary, Arch. Rational Mech. Anal., {\bf 92}, 1986, pp.205-245.
\item[{[3]}] Friedman, A., The Stefan problem in several space variables,  Trans. Amer. Math. Soc.  {\bf 133}, 1968, pp. 51-87.
\item[{[4]}] Friedman, A., Partial differential equations of parabolic type, Prentice-Hall, Englewood Cliffs, N. J., 1964.
\item[{[5]}] Friedman, A., Kinderlehrer, D., A., one phese Stefan problem, Indiana Univ. Math. J. {\bf 24}, 1975, pp. 1005-1035.
\item[{[6]}] Kamenomostskaya, S. L., On Stefan problem, Math. Sb., {\bf 53(95)}, No.4, 1961, pp 489-514 (in Russian).
\item[{[7]}] Lam\'e, G, Clapeyron, B, P, M\'emoire sur la solidification par refroidissement d'un globe solide, Ann. Chem. Phys., {\bf 47}, 1831, pp.250-256.
\item[{[8]}] Reutskiy, S., Y., The method of approximate fundamental solutions (MAFS) for Stefan problems for the sphere. Appl. Math. Comput., {\bf 227}, 2014, pp. 648-655.
\item[{[9]}] Salsa. S., Two-phase Stefan problems. Recent results and open questions. Milan J. Math. {\bf 80}, 2012, pp 267-281. 
\item[{[10]}] Yi, F., Global classical solution of Muskat free boundary problem. J. Math. Anal. Appl.  {\bf 288}, 2003, pp 442-461.
\end{itemize}
\end{document}